\documentstyle[leqno]{article}

\newtheorem{th}{Theorem}[section]
\newtheorem{prop}[th]{Proposition}

\newcounter{defin}[section]
\renewcommand{\thedefin}{\thesection.\arabic{defin}}
\newcounter{ex}[section]

\newcounter{rem}[section]
\renewcommand{\therem}{\thesection.\arabic{rem}}
\title{The Generalized Metrical Multi-Time Lagrange \\
Space  of Relativistic Geometrical Optics}
\date{}
\author{Mircea Neagu}
\begin{document}
\maketitle
\begin{abstract}
Section 1 contains some physical and geometrical aspects that motivates us to
study the generalized  metrical multi-time Lagrange space of Relativistic
Geometrical Optics, denoted by $RGOGML^n_p$.
Section 2  developes the geometry  of this space, in the
sense of d-connections, d-torsions and d-curvatures. 
The  Einstein equations of gravitational potentials of this  generalized metrical
multi-time Lagrange space are studied in Section 3. The conservation laws of
the stress-energy d-tensor  of $RGOGML^n_p$  are also described.
The electromagnetic d-tensors are introduced in Section 4, and corresponding
Maxwell equations  are derived.
\end{abstract}
{\bf Mathematics Subject Classification (2000):} 53B40, 53C60, 53C80.\\
{\bf Key words:} 1-jet fibre bundle, nonlinear connection, Cartan canonical
connection, Einstein equations, Maxwell equations.

\section{Geometrical and physical aspects}

\hspace{5mm} Let us consider the generalized Lagrange space $GL^n=(M,
g_{ij}(x,y)),\;n=\dim  M$,
whose fundamental tensor field $g_{ij}(x,y)$ on $TM$ is of the form
\begin{equation}\label{fdt}
g_{ij}(x,y)=\varphi_{ij}(x)+\left[1-{1\over n^2(x,y)}\right]y_iy_j,
\end{equation}
where $\varphi_{ij}(x)$ is a semi-Riemannian metric tensor on $M$, $y_i=\varphi
_{ij}(x)y^j$ and $n(x,y)$ is a smooth function on $TM$, which  is  called  the
{\it refractive index} function. This generalized Lagrange space is known today as the
generalized Lagrange space of {\it relativistic  geometrical optics} \cite{9}.

In  order to  explain the above  physical terminology, let us analyse the restriction
of the geometry of the previous space to a cross section,
\begin{equation}
S_V:M\to TM,\quad x\to(x,y^i=V^i(x)).
\end{equation}
Thus, starting with a given cross section $S_V$, we remark that
the restriction of the fundamental d-tensor
$g_{ij}(x,y)$ of $GL^n$ to the submanifold $S_V(M)$ is given by
\begin{equation}\label{synge}
g_{ij}(x,V(x))=\varphi_{ij}(x)+\left[1-{1\over n^2(x,V(x))}\right]V_iV_j,
\end{equation}
where $V_i=\varphi_{ij}(x)V^j$. We underline that the metric \ref{synge} was
introduced by Synge \cite{20} and  was used by him in the study of the propagation of the
electromagnetic waves in a medium with the index of refraction $n(x,V(x))$,
$V(x)$ being the velocity of medium. In Synge's terminology, a {\it medium} is
represented by a triad  {\bf M}$=(M,V(x),n(x,V(x))$. More deeply, he use the
following
\medskip\\
\addtocounter{defin}{1}
{\bf Definition \thedefin} A triad  {\bf M}$=(M,V(x),n(x,V(x))$ is called
a {\it dispersive medium}. If $n(x,y)$ does not depend  of $y$, then  {\bf M}
is called a {\it  non-dispersive medium}.\medskip

Consequently,  the geometry induced by the  Synge's metric tensor \ref{synge}
gives a mathematical model  for  relativistic optics.

In conclusion, the study  of the generalized  Lagrange space $GL^n$, whose
metrical d-tensor is  of the form  \ref{fdt}, was imposed. In this direction,
an important  class  of generalized Lagrange spaces of relativistic geometrical
optics,  having  the  refractive index in the form,
\begin{equation}\label{pc}
\displaystyle{{1\over n^2}=1-{\alpha\over  c^2},\;\alpha\in R^*_+},
\end{equation}
$c$ being the light velocity, was studied by Miron and Kawaguchi \cite{11}.
For these spaces, some post-Newtonian  estimations were investigated by Asanov
and  Kawaguchi \cite{1}, \cite{7}. In a different version, some post-Newtonian
estimations was presented also  by Roxburgh \cite{18}.

The geometry of the generalized Lagrange spaces of relativistic geometrical
optics is now completely done by Miron ,Anastasiei and Kawaguchi \cite{9},
\cite{11}. Their geometrical development relies the using of an {\it "a priori"}
fixed nonlinear connection on $TM$, whose
components are
\begin{equation}\label{mnc}
N^i_j(x,y)=\gamma^i_{jk}(x)y^k,
\end{equation}
where $\gamma^i_{jk}(x)$ are the Christoffel symbols for the semi-Riemannian
metric $\varphi_{ij}(x)$.

The  using of the nonlinear connection \ref{mnc} in the study of the generalized Lagrange
space of relativistic geometrical optics $GL^n$ is motivated in various ways.
In this direction, we present only two  geometrical and physical aspects. For
more details, see \cite{9}, \cite{11}.

Firstly, it is very  important that the  autoparallel curves of the nonlinear
connection $N^i_j(x,y)$ of the generalized Lagrange space of relativistic
geometrical optics $GL^n$ coincide to  the geodesics of the  Riemannian space
$R^n=(M,\varphi_{ij})$. In other words, the space $GL^n$ verifies the first
EPS (Ehlers, Pirani, Schild) condition from the constructive-axiomatic
formulation of General Relativity \cite{4}.

Secondly, it is remarkable that, in the particular case \ref{pc}, the
{\it absolute energy} Lagrangian of $GL^n$,
\begin{equation}
{\cal E}:TM\to R,\quad {\cal E}(x,y)=g_{ij}(x,y)y^iy^j,
\end{equation}
is a {\it  regular Lagrangian} and then its canonical nonlinear connection
\cite{11} is exactly that given  by \ref{mnc}.

In conclusion, one can assert that a generalized Lagrange space $GL^n$, which
verifies the above axiomatic assumptions,  becomes  a convenient mathematical
model for the relativistic geometrical optics.

In this paper, we try to extend the previous geometrical and physical theories,
from the  tangent bundle  $TM$ to the more general jet fibre bundle of order
one $J^1(T,M)$ coordinated by $(t^\alpha,x^i,x^i_\alpha)$, where
$T$  is a smooth, real, $p$-dimensional manifold coordinated by $t=(t^\alpha)
_{\alpha=\overline{1,p}}$, whose physical meaning is that of {\it
"multidimensional time"}, while $x^i_\alpha$ have the physical meaning of {\it
partial directions}. In this sense, we recall that the jet fibre
bundle of order one $J^1(T,M)$ is a basic object in the study of classical
and quantum field theories.

A natural geometry  of physical fields induced by a Kronecker $h$-regular
vertical metrical multi-time d-tensor $G^{(\alpha)(\beta)}_{(i)(j)}(t^\gamma,x^k,
x^k_\gamma)$ on the total space  of the 1-jet vector bundle $J^1(T,M)\to
T\times M$, where $h=(h_{\alpha\beta}(t^\gamma))$ is a semi-Riemannian metric
on the temporal manifold $T$,  was created by Neagu \cite{12}.

The fundamental geometrical concept used there is that of {\it
generalized metrical multi-time Lagrange space}. This geometrical concept
with physical meaning is represented by a pair $GML^n_p=(J^1(T,M),
G^{(\alpha)(\beta)}_{(i)(j)})$ consisting of the 1-jet space and a
{\it Kronecker $h$-regular} vertical multi-time metrical d-tensor
$G^{(\alpha)\beta)}_{(i)(j)}$ on $J^1(T,M)$, that  is, it decomposes in
\begin{equation}
G^{(\alpha)(\beta)}_{(i)(j)}(t^\gamma,x^k,x^k_\gamma)=h^{\alpha\beta}(t^
\gamma)g_{ij}(t^\gamma,x^k,x^k_\gamma),
\end{equation}
where $g_{ij}(t^\gamma,x^k,x^k_\gamma)$ is a d-tensor on $J^1(T,M)$,
symmetric, of rank $n$ and having a constant signature. The d-tensor
$g_{ij}(t^\gamma,x^k,x^k_\gamma)$ is called the
{\it spatial metrical d-tensor of $GML^n_p$}.

Following the general physical and geometrical development from \cite{12},
the aim of this paper is to study the particular generalized metrical
multi-time Lagrange space $RGOGML^n_p$, whose spatial metrical d-tensor is of
the form
\begin{equation}
g_{ij}(t^\gamma,x^k,x^k_\gamma)=\varphi_{ij}(x^k)+A_i(t^\gamma,x^k,x^k_\gamma)
A_j(t^\gamma,x^k,x^k_\gamma)
\end{equation}
where $\varphi_{ij}(x^k)$  is  a semi-Riemannian metric on the spatial
manifold  $M$ and $A_i(t^\gamma,x^k,x^k_\gamma)$ represent the components
of a d-tensor $A$ on $J^1(T,M)$, whose physical meaning  is that of {\it
refractive index} d-tensor of the  {\it medium} {\bf M}=$J^1(T,M)$.
\medskip\\
\addtocounter{rem}{1}
{\bf Remarks \therem} i) To motivate the terminology used above, let us consider
the particular case when the  temporal manifold identifies with  the usual
time axis, represented   by the  set of real numbers $R$. In that case,
setting
\begin{equation}
A_i(t,x^k,y^k)=\sqrt{1-{1\over n^2(t,x^k,y^k)}}\cdot y_i,
\end{equation}
where $y_i=\varphi_{im}(x^k)y^m$ and $n:J^1(R,M)\equiv R\times TM\to[1,\infty)$
is a smooth function, we can regard $n$ like a  {\it dynamic refractive index}
function (i. e. the refractive index modifies in time). Consequently, this
particular space represents a natural dynamical generalization of the classical
generalized Lagrange space  of relativistic geometrical optics from \cite{9},
\cite{11}.

ii) The  inverse of the spatial metrical d-tensor $g_{ij}$ of $RGOGML^n_p$ is given
by  the following d-tensor field,
\begin{equation}\label{inv}
g^{ij}=\varphi^{ij}-{1\over 1+A_0}A^iA^j,
\end{equation}
where $A^i=\varphi^{im}A_m$ and  $A_0=A^mA_m$.
\medskip

To develope the geometry of this generalized metrical multi-time
Lagrange space, we need  a nonlinear connection $\Gamma=(M^{(i)}_{(\alpha)\beta},
N^{(i)}_{(\alpha)j})$ on $J^1(T,M)$ \cite{14}. In this direction, we fix {\it
"a priori"} the nonlinear  connection $\Gamma$ defined by the temporal components
\begin{equation}
M^{(i)}_{(\alpha)\beta}=-H^\mu_{\alpha\beta}x^i_\mu
\end{equation}
and the spatial components
\begin{equation}\label{nnc}
N^{(i)}_{(\alpha)j}=\gamma^i_{jm}x^m_\alpha,
\end{equation}
where $H^\alpha_{\beta\gamma}$ (resp. $\gamma^i_{jk}$) are the Christoffel
symbols  of the semi-Riemannian metric $h_{\alpha\beta}$ (resp. $\varphi_{ij}$).
\medskip\\
\addtocounter{rem}{1}
{\bf Remarks \therem} i) Our given {\it  a priori} nonlinear connection $\Gamma$
on $J^1(T,M)$ naturally  generalizes that  used   by  Miron, Anastasiei  and
Kawaguchi.

ii) The spatial components $N^{(i)}_{(\alpha)j}$ of the fixed nonlinear
connection $\Gamma$ are  {\it without torsion} \cite{12}.

iii) The previous nonlinear connection $\Gamma$  is dependent
only the vertical fundamental metrical d-tensor $G^{(\alpha)(\beta)}_{(i)(j)}$
of $RGOGML^n_p$. This fact emphasize the {\it metrical character} of the geometry
attached to this space, i. e. , all geometrical  objects are directly arised from
$G^{(\alpha)(\beta)}_{(i)(j)}$.

iv) Using the relation between sprays and the components of a nonlinear
connection and the definition of harmonic maps attached to a given multi-time
dependent spray on $J^1(T,M)$ (for more  details, see \cite{14}), we easily
deduce that the harmonic maps of the nonlinear connection $\Gamma$ of $RGOGML^n_p$
are exactly the harmonic maps between the semi-Riemannian spaces $(T,h)$ and
$(M,\varphi)$ \cite{3}.
\medskip

In conclusion, we can assert that the generalized metrical multi-time Lagrange
space $RGOGML^n_p$, which verifies the previous assumptions, represents a
convenient geometrical model for relativistic geometrical optics, in a general
setting.\medskip\\
{\bf Open problem.} At the  end of this section, we should like to point out
that, in the particular case $A_i=A_i(t^\gamma,x^k)$ (i. e. the refractive
index d-tensor $A$ does not depend by partial directions
$x^i_\alpha$), the {\it absolute energy} Lagrangian function \cite{12},
\begin{equation}
{\cal E}:J^1(T,M)\to R,\quad {\cal E}(t^\gamma,x^k,x^k_\gamma)=h^{\alpha\beta}
(t^\gamma)g_{ij}(t^\gamma,x^k)x^i_\alpha x^j_\beta,
\end{equation}
is a {\it Kronecker $h$-regular Lagrangian} \cite{15}. Consequently, it naturally
induces a canonical spatial nonlinear connection on $J^1(T,M)$, whose
components are given by \cite{15}
\begin{equation}\label{bnnc}
\bar N^{(i)}_{(\alpha)j}=\Gamma^i_{jm}x^m_\alpha+{g^{im}\over 2}{\partial
g_{mj}\over\partial t^\alpha},
\end{equation}
where
$$
\displaystyle{\Gamma^l_{jk}={g^{li}\over 2}\left({\partial g_{ij}\over\partial
x^k}+{\partial g_{ik}\over\partial x^j}-{\partial g_{jk}\over\partial x^i}
\right)}
$$
are the {\it generalized Christoffel symbols of the "multi-time" dependent
spatial metric $g_{ij}$}. In conclusion, it is natural to arise the following
question:\medskip

$-$Considering  an {\it  isotropic medium} {\bf M}=$J^1(T,M)$ (i. e. its refractive
index  d-tensor $A$ does  not depend by partial directions $x^i_\alpha$), what
is the difference, from a physical point of view, between the  using  of
one or  another one of the spatial nonlinear connections expressed by
\ref{nnc} and \ref{bnnc} ?

\section{Cartan canonical connection}

\setcounter{equation}{0}
\hspace{5mm} In this section, we will apply the general geometrical development
of the generalized metrical multi-time Lagrange spaces \cite{12}, to the
particular space of relativistic geometrical optic
$RGOGML^n_p=(J^1(T,M),G^{(\alpha)(\beta)}_{(i)(j)})$, having the vertical
fundamental d-tensor,
\begin{equation}\label{mdt}
G^{(\alpha)(\beta)}_{(i)(j)}(t^\gamma,x^k,x^k_\gamma)=h^{\alpha\beta}
(t^\gamma)\left[\varphi_{ij}(x^k)+A_i(t^\gamma,x^k,x^k_\gamma)A_j(t^\gamma,
x^k,x^k_\gamma)\right],
\end{equation}
and being endowed with the nonlinear  connection $\Gamma=(M^{(i)}_{(\alpha)\beta},
N^{(i)}_{(\alpha)j})$, where
\begin{equation}\label{nc}
M^{(i)}_{(\alpha)\beta}=-H^\mu_{\alpha\beta}x^i_\mu,
\qquad
N^{(i)}_{(\alpha)j}=\gamma^i_{jm}x^m_\alpha.
\end{equation}

Let $\displaystyle{\left\{{\delta\over\delta
t^\alpha}, {\delta\over\delta x^i}, {\partial\over\partial x^i_\alpha}\right\}
\subset{\cal X}(J^1(T,M))}$ and $\{dt^\alpha, dx^i, \delta x^i_\alpha\}\subset{\cal
X}^*(J^1(T,M))$  be the adapted bases of the nonlinear connection $\Gamma$,
where \cite{14}
\begin{equation}
\left\{\begin{array}{l}\medskip
\displaystyle{{\delta\over\delta t^\alpha}={\partial\over\partial t^\alpha}-
M^{(j)}_{(\beta)\alpha}{\partial\over\partial x^j_\beta}}\\\medskip
\displaystyle{{\delta\over\delta x^i}={\partial\over\partial x^i}-
N^{(j)}_{(\beta)i}{\partial\over\partial x^j_\beta}}\\
\delta x^i_\alpha=dx^i_\alpha+M^{(i)}_{(\alpha)\beta}dt^\beta+N^{(i)}_{(\alpha)
j}dx^j.
\end{array}\right.
\end{equation}

Following the paper \cite{12}, by direct computations, we can  determine the
{\it Cartan canonical connection} of $RGOGML^n_p$, together with its torsion
and curvature local d-tensors.

In order to describe  these  geometrical entities of $RGOGML^n_p$,
let us consider $B\Gamma=(H^\gamma_{\alpha\beta},0,\gamma^i_{jk},0)$, the
Berwald $h$-normal $\Gamma$-linear connection attached to the semi-Riemannian
metrics $h_{\alpha\beta}$ and $\varphi_{ij}$ and $"_{//\alpha}"$, $"_{\Vert_i}"$,
$"\Vert^{(\alpha)}_{(i)}"$, its local covariant derivatives \cite{13}.
It is easy to deduce that the Berwald connection $B\Gamma$ of $RGOGML^n_p$
has the following metrical properties:
\begin{equation}
\left\{\begin{array}{lll}\medskip
h_{\alpha\beta//\gamma}=0,\quad h_{\alpha\beta\Vert k}=0,\quad
h_{\alpha\beta}\Vert^{(\gamma)}_{(k)}=0,\\\medskip
\varphi_{ij//\gamma}=0,\quad \varphi_{ij\Vert k}=0,\quad
\varphi_{ij}\Vert^{(\gamma)}_{(k)}=0.
\end{array}\right.
\end{equation}

In this context, using the general expressions which give the components of
the  Cartan canonical connection of a generalized metrical multi-time
Lagrange space \cite{12}, by a direct calculation, we obtain

\pagebreak

\begin{th}
The Cartan canonical connection
$
C\Gamma=(H^\gamma_{\alpha\beta}, G^i_{j\gamma}, L^i_{jk}, C^{i(\gamma)}_{j(k)})
$
of $RGOGML^n_p$ has the adapted coefficients
\begin{equation}\label{cc}
\begin{array}{l}\medskip
H^\gamma_{\alpha\beta}=H^\gamma_{\alpha\beta},\\\medskip
\displaystyle{G^i_{j\gamma}={1\over 2}\left[(A^iA_j)_{//\gamma}-{A^iA^m\over
1+A_0}(A_mA_j)_{//\gamma}\right],}\\\medskip
\displaystyle{L^i_{jk}=(\gamma_{jkm}+A_{jkm})\left[\varphi^{im}-{A^iA^m\over
1+A_0}\right]},\\
\displaystyle{C^{i(\gamma)}_{j(k)}=C^\gamma_{jkm}\left[\varphi^{im}-{A^iA^m
\over 1+A_0}\right]},
\end{array}
\end{equation}
where $A^i=\varphi^{im}A_m$, $A_0=A^mA_m$ and
\begin{equation}
\begin{array}{l}\medskip
\displaystyle{\gamma_{ijm}={1\over 2}\left[{\partial\varphi_{im}\over\partial
x^j}+{\partial\varphi_{jm}\over\partial x^i}-{\partial\varphi_{ij}\over
\partial x^m}\right]},\\\medskip
\displaystyle{A_{ijm}={1\over 2}\left[{\delta(A_iA_m)\over\delta x^j}+
{\delta(A_jA_m)\over\delta x^i}-{\delta(A_iA_j)\over\delta x^m}\right]},\\
\displaystyle{C^\gamma_{ijm}={1\over 2}\left[{\partial(A_iA_m)\over\partial
x^j_\gamma}+{\partial(A_jA_m)\over\partial x^i_\gamma}-{\partial(A_iA_j)\over
\partial x^m_\gamma}\right]}.
\end{array}
\end{equation}
\end{th}
\addtocounter{rem}{1}
{\bf Remark \therem} Using the  notations
\begin{equation}
\begin{array}{ll}\medskip
A^i_{jk}=\varphi^{im}A_{jkm},&\displaystyle{\Lambda^i_{jk}=A^i_{jk}-{(\gamma_
{jk0}+A_{jk0})A^i\over 1+A_0}},\\
\stackrel{1}{C}^{\gamma i}_{jk}=\varphi^{im}C^\gamma_{jkm},&\displaystyle{
\stackrel{0}{C}^{\gamma i}_{jk}=-{C^\gamma_{jk0}A^i\over 1+A_0}},
\end{array}
\end{equation}
where $D^{\ldots}_{\ldots 0}=D^{\ldots}_{\ldots m}A^m$, we can rewrite  the
components $L^i_{jk}$ and $C^{i(\gamma)}_{j(k)}$ of $C\Gamma$ in the
following simple form:
\begin{equation}\label{ncc}
\begin{array}{l}\medskip
L^i_{jk}=\gamma^i_{jk}+\Lambda^i_{jk},\\
C^{i(\gamma)}_{j(k)}=\stackrel{1}{C}^{\gamma i}_{jk}+\stackrel{0}{C}^{\gamma
i}_{jk}.
\end{array}
\end{equation}

\begin{th}
The torsion {\em\bf T} of the Cartan canonical connection of $RGOGML^n_p$
is determined by seven effective local d-tensors, namely,
\begin{equation}
\begin{array}{l}\medskip
T^m_{\alpha j}=-G^m_{j\alpha},\quad P^{(m)\;\;(\beta)}_{(\mu)\alpha(j)}=-
\delta^\beta_\mu G^m_{j\alpha},\quad P^{(m)\;(\beta)}_{(\mu)i(j)}=-\delta^
\beta_\mu\Lambda^m_{ij},
\\\medskip
P^{m(\beta)}_{i(j)}=\stackrel{1}{C}^{\beta m}_{ij}+\stackrel{0}{C}^{\beta m}
_{ij},\quad
S^{(m)(\alpha)(\beta)}_{(\mu)(i)(j)}=\delta^\alpha_\mu C^{m(\beta)}_{i(j)}-
\delta^\beta_\mu C^{m(\alpha)}_{i(j)},\\
R^{(m)}_{(\mu)\alpha\beta}=-H^\gamma_{\mu\alpha\beta}x^m_\gamma,
\quad R^{(m)}_{(\mu)\alpha j}=0,\quad
R^{(m)}_{(\mu)ij}=r^m_{ijk}x^k_\mu,
\end{array}
\end{equation}
where $H^\gamma_{\mu\alpha\beta}$ (resp. $r^m_{ijk}$) are the local curvature
tensors of the semi-Riemannian metric $h_{\alpha\beta}$ (resp. $\varphi_{ij}$).
\end{th}

The general expressions of  the local curvature d-tensors attached to the
Cartan canonical connection of a generalized  metrical multi-time Lagrange
space \cite{12}, applied to $RGOGML^n_p$, imply

\begin{th}
The curvature {\em\bf R} of the Cartan canonical connection of $RGOGML^n_p$
is determined by seven effective local d-tensors, expressed by,

\begin{equation}\label{curv}
\begin{array}{l}\medskip
\displaystyle{H^\alpha_{\eta\beta\gamma}={\partial H^\alpha_{\eta\beta}\over
\partial t^\gamma}-{\partial H^\alpha_{\eta\gamma}\over\partial t^\beta}+
H^\mu_{\eta\beta}H^\alpha_{\mu\gamma}-H^\mu_{\eta\gamma}H^\alpha_{\mu\beta},}
\\\medskip
\displaystyle{R^l_{i\beta\gamma}={\delta G^l_{i\beta}\over\delta t^\gamma}-
{\delta G^l_{i\gamma}\over\delta t^\beta}+G^m_{i\beta}G^l_{m\gamma}-
G^m_{i\gamma}G^l_{m\beta}-H^\eta_{\mu\beta\gamma}C^{l(\mu)}_{i(m)}x^m_\eta,}
\\\medskip
\displaystyle{R^l_{i\beta k}={\delta G^l_{i\beta}\over\delta x^k}-
{\delta L^l_{ik}\over\delta t^\beta}+G^m_{i\beta}L^l_{mk}-
L^m_{ik}G^l_{m\beta},}
\\\medskip
\displaystyle{R^l_{ijk}=r^l_{ijk}+\rho^l_{ijk},}
\\\medskip
\displaystyle{P^{l\;\;(\gamma)}_{i\beta(k)}={\partial G^l_{i\beta}\over\partial
x^k_\gamma}-C^{l(\gamma)}_{i(k)/\beta}-C^{l(\gamma)}_{i(m)}G^m_{k\beta},}
\\\medskip
\displaystyle{P^{l\;(\gamma)}_{ij(k)}={\partial\Lambda^l_{ij}\over\partial
x^k_\gamma}-\Lambda^m_{jk}\stackrel{1}{C}^{\gamma l}_{im}-\Lambda^m_{jk}
\stackrel{0}{C}^{\gamma l}_{im}-\stackrel{1}{C}^{\gamma l}_{ik\vert j}-
\stackrel{0}{C}^{\gamma l}_{ik\vert j},}
\\\medskip
S^{l(\beta)(\gamma)}_{i(j)(k)}=\stackrel{1}{S}^{l(\beta)(\gamma)}_{i(j)(k)}-
\stackrel{0}{S}^{l(\beta)(\gamma)}_{i(j)(k)},
\end{array}
\end{equation}
where  $H^\eta_{\alpha\beta\gamma}$ and $r^l_{ijk}$ are the curvature tensors
of the semi-Riemannian metrics $h_{\alpha\beta}$ and $\varphi_{ij}$, the
operators $"_{/\beta}"$, $"_{\vert j}"$ represent the  local covariant
derivatives of $C\Gamma$, and
\begin{equation}
\begin{array}{l}\medskip
\displaystyle{\stackrel{1}{S}^{l(\beta)(\gamma)}_{i(j)(k)}={\partial\stackrel{1}
{C}^{\beta  l}_{ij}\over\partial x^k_\gamma}-{\partial\stackrel{1}{C}^{
\gamma l}_{ik}\over\partial x^j_\beta}+\stackrel{1}{C}^{\beta m}_{ij}
C^{l(\gamma)}_{m(k)}-\stackrel{1}{C}^{\beta l}_{mj}C^{m(\gamma)}_{i(k)},}\\\medskip
\displaystyle{\stackrel{0}{S}^{l(\beta)(\gamma)}_{i(j)(k)}={\partial\stackrel{0}
{C}^{\beta l}_{ij}\over\partial x^k_\gamma}-{\partial\stackrel{0}{C}^{
\gamma l}_{ik}\over\partial x^j_\beta}+\stackrel{0}{C}^{\beta m}_{ij}
C^{l(\gamma)}_{m(k)}-\stackrel{0}{C}^{\beta l}_{mj}C^{m(\gamma)}_{i(k)},}\\
\rho^l_{ijk}=\Lambda^l_{ij\vert k}-\Lambda^l_{ik\vert j}+\Lambda^m_{ij}
\Lambda^l_{mk}-\Lambda^m_{ik}\Lambda^l_{mj}+C^{l(\mu)}_{i(m)}r^m_{sjk}x^s_\mu.
\end{array}
\end{equation}
\end{th}

\section{Einstein equations of gravitational field}

\setcounter{equation}{0}

\hspace{5mm} Concerning the gravitational theory on $RGOGML^n_p$,
we point out that the vertical metrical d-tensor \ref{mdt} and its fixed
nonlinear connection \ref{nc} induce a natural {\it gravitational
$h$-potential} on the 1-jet space $J^1(T,M)$ (i.  e. a Sasakian-like metric),
which is expressed by \cite{12}
\begin{equation}
G=h_{\alpha\beta}dt^\alpha\otimes dt^\beta+g_{ij}dx^i\otimes
dx^j+h^{\alpha\beta}g_{ij}\delta x^i_\alpha\otimes\delta x^j
_\beta,
\end{equation}
where $g_{ij}=\varphi_{ij}+A_iA_j$.
Let $C\Gamma=(H^\gamma_{\alpha\beta},G^k_{j\gamma},L^i_{jk},C^{i(\gamma)}_
{j(k)})$ be the Cartan canonical connection of $RGOGML^n_p$.

We postulate that the Einstein equations which govern the gravitational
$h$-potential $G$ of $RGOGML^n_p$ are the Einstein equations attached to
the Cartan canonical
connection and the adapted metric $G$ on $J^1(T,M)$, that is,
\begin{equation}
Ric(C\Gamma)-{Sc(C\Gamma)\over 2}G={\cal K}{\cal T},
\end{equation}
where $Ric(C\Gamma)$ represents the Ricci d-tensor of the Cartan connection,
$Sc(C\Gamma)$ is its scalar curvature, ${\cal K}$ is the Einstein constant
and ${\cal T}$ is an intrinsec distinguished tensor of matter which is called  the {\it
stress-energy} d-tensor.

In the adapted basis $(X_A)=\displaystyle{\left({\delta\over\delta t^\alpha},
{\delta\over\delta x^i},{\partial\over\partial x^i_\alpha}\right)}$ attached
to $\Gamma$, the curvature d-tensor {\bf R} of the Cartan connection is
expressed locally by {\bf R}$(X_C,X_B)X_A=R^D_{ABC}X_D$. Hence, it follows
that we have $R_{AB}=Ric(C\Gamma)(X_A,X_B)=R^D_{ABD}$ and $Sc(C\Gamma)=G^{AB}
R_{AB}$, where
\begin{equation}
G^{AB}=\left\{\begin{array}{ll}\medskip
h_{\alpha\beta},&\mbox{for}\;\;A=\alpha,\;B=\beta\\\medskip
g^{ij},&\mbox{for}\;\;A=i,\;B=j\\\medskip
h_{\alpha\beta}g^{ij},&\mbox{for}\;\;A={(i)\atop(\alpha)},\;B={(j)\atop(\beta)}\\
0,&\mbox{otherwise},
\end{array}\right.
\end{equation}
the tensor field $g^{ij}$ being expressed by \ref{inv}.

Taking into account the expressions \ref{curv} of the local curvature
d-tensors of the Cartan connection  of $RGOGML^n_p$, we obtain without difficulties
\begin{th}
The Ricci  d-tensor $Ric(C\Gamma)$ of $RGOGML^n_p$ is determined by seven effective
local d-tensors expressed, in adapted basis, by:
\begin{equation}
\hspace*{5mm}
\begin{array}{l}\medskip
H_{\alpha\beta}=H^\mu_{\alpha\beta\mu},
\quad
R_{i\beta}=R^m_{i\beta m},
\quad
R_{ij}=r_{ij}+\rho_{ij},
\quad
P^{(\alpha)}_{(i)\beta}=P^{m\;(\alpha)}_{i\beta(m)},
\\\medskip
P^{\;(\alpha)}_{i(j)}=-P^{m\;(\alpha)}_{im(j)},
\quad
P^{(\alpha)}_{(i)j}=-P^{\;(\alpha)}_{i(j)},
\quad
S^{(\beta)(\gamma)}_{(j)(k)}=\stackrel{1}{S}^{(\beta)(\gamma)}_{(j)(k)}-
\stackrel{0}{S}^{(\beta)(\gamma)}_{(j)(k)},
\end{array}
\end{equation}
where $H_{\alpha\beta}$ (resp. $r_{ij}$) are the local Ricci tensors of the
semi-Riemannian metric $h_{\alpha\beta}$ (resp. $\varphi_{ij}$),
$\rho_{ij}=\rho^m_{ijm}$, $\stackrel{1}{S}^{(\beta)(\gamma)}_{(j)(k)}=
\stackrel{1}{S}^{m(\gamma)(\beta)}_{i(j)(m)}$ and
$\stackrel{0}{S}^{(\beta)(\gamma)}_{(j)(k)}=\stackrel{0}{S}^{m(\gamma)(\beta)}
_{i(j)(m)}$.
\end{th}

Let us denote $H=h^{\alpha\beta}H_{\alpha\beta}$, $R=g^{ij}R_{ij}$ and
$S=h_{\alpha\beta}g^{ij}S^{(\alpha)(\beta)}_{(i)(j)}$. In this context,
by a simple  calculation, it follows
\begin{th}
The scalar curvature of the Cartan connection $C\Gamma$ of $RGOGML^n_p$ has
the formula
\begin{equation}
Sc(C\Gamma)=H+R+S,
\end{equation}
the terms $H,\;R$ and $S$ being determined by the relations:
\begin{equation}
\begin{array}{l}\medskip
H=h^{\alpha\beta}H_{\alpha\beta},\\\medskip
\displaystyle{R=r+\rho-{r_{00}+\rho_{00}\over  1+A_0},}\\
S=\stackrel{1}{S}-\stackrel{1}{S^\prime}+\stackrel{0}{S}-\stackrel{0}{S
^\prime},
\end{array}
\end{equation}
where $H$ (resp $r$) is the scalar curvature of the semi-Riemannian
metric $h_{\alpha\beta}$ (resp. $\varphi_{ij}$) and
\begin{equation}
\begin{array}{l}\medskip
\rho=\varphi^{rs}\rho_{rs},\quad r_{00}=r_{ms}A^mA^s,\quad
\rho_{00}=\rho_{ms}A^mA^s,\\\medskip
\displaystyle{\stackrel{1}{S}=h_{\mu\nu}\varphi^{rs}\stackrel{1}{S}^{(\mu)(\nu)}
_{(r)(s)},\quad \stackrel{1}{S^\prime}={h_{\mu\nu}A^rA^s\stackrel{1}{S}
^{(\mu)(\nu)}_{(r)(s)}\over 1+A_0}},\\
\displaystyle{\stackrel{0}{S}=h_{\mu\nu}\varphi^{rs}\stackrel{0}{S}^{(\mu)(\nu)}
_{(r)(s)},\quad \stackrel{0}{S^\prime}={h_{\mu\nu}A^rA^s\stackrel{0}{S}
^{(\mu)(\nu)}_{(r)(s)}\over 1+A_0}}.
\end{array}
\end{equation}
\end{th}

Following the gravitational field theory exposition on a generalized metrical
multi-time Lagrange space $GML^n_p$, $\dim  M=n,\;\dim T=p$, from the paper
\cite{12}, by  local computations,
we can give
\begin{th}
If $p>2$ and $n>2$, the Einstein equations which govern the gravitational
$h$-potential $G$ of $RGOGML^n_p$ have  the  local form
$$
\left\{\begin{array}{l}\medskip
\displaystyle{H_{\alpha\beta}-{H\over 2}h_{\alpha\beta}={\cal K}\tilde{\cal T}_
{\alpha\beta}}\\\medskip
\displaystyle{r_{ij}-{r\over 2}\varphi_{ij}+\theta_{ij}={\cal K}\tilde{\cal T}_{ij}}\\
\displaystyle{S^{(\alpha)(\beta)}_{(i)(j)}+{\stackrel{1}{S}-\stackrel{1}
{S^\prime} \over 2}h^{\alpha\beta}g_{ij}+\stackrel{00}{S}^{(\alpha)(\beta)}
_{(i)(j)}={\cal K}\tilde{\cal T}^{(\alpha)(\beta)}_{(i)(j)}},
\end{array}\right.\leqno{(E^\prime_1)}
$$
$$
\left\{\begin{array}{lll}\medskip
0={\cal T}_{\alpha i},&R_{i\alpha}={\cal K}{\cal T}_{i\alpha},&
P^{(\alpha)}_{(i)\beta}={\cal K}{\cal T}^{(\alpha)}_{(i)\beta}\\
0={\cal T}^{\;(\beta)}_{\alpha(i)},&
P^{\;(\alpha)}_{i(j)}={\cal K}{\cal T}^{\;(\alpha)}_{i(j)},&
P^{(\alpha)}_{(i)j}={\cal K}{\cal T}^{(\alpha)}_{(i)j},
\end{array}\right.\leqno{(E_2)}
$$
where $\tilde{\cal T}_{\alpha\beta}$, $\tilde{\cal T}_{ij}$ and $\tilde{\cal
T}^{(\alpha)(\beta)}_{(i)(j)}$  represent the components  of a new
stress-energy d-tensor $\tilde{\cal T}$, defined by the relations
\begin{equation}\label{*}
\left\{\begin{array}{l}\medskip
\displaystyle{\tilde{\cal T}_{\alpha\beta}={\cal T}_{\alpha\beta}+{R+S\over
2{\cal K}}h_{\alpha\beta}}\\\medskip
\displaystyle{\tilde{\cal T}_{ij}={\cal T}_{ij}+{H+S\over 2{\cal K}}
g_{ij}}\\
\displaystyle{\tilde{\cal T}^{(\alpha)(\beta)}_{(i)(j)}={\cal T}^{(\alpha)
(\beta)}_{(i)(j)}+{H+R\over 2{\cal K}}h^{\alpha\beta}g_{ij}},
\end{array}\right.
\end{equation}
and
\begin{equation}
\begin{array}{l}\medskip
\displaystyle{\theta_{ij}=\rho_{ij}-{1\over 2}\left(\rho-{r_{00}+\rho_{00}\over
1+A_0}\right)g_{ij},}\\
\displaystyle{\stackrel{00}{S}^{(\alpha)(\beta)}_{(i)(j)}=\stackrel{0}{S}^
{(\alpha)(\beta)}_{(i)(j)}-{\stackrel{0}{S}-\stackrel{0}{S^\prime}\over 2}
h^{\alpha\beta}g_{ij}.}
\end{array}
\end{equation}
\end{th}
\addtocounter{rem}{1}
{\bf Remark \therem} Note that, in order to have the compatibility of the Einstein equations, it is
necessary that the certain  adapted local components of the stress-energy
d-tensor vanish {\it "a priori"}.
\medskip

From physical point of view, the
stress-energy d-tensor ${\cal T}$ must  verify the local {\it conservation
laws} ${\cal T}^B_{A\vert B}=0,\;\forall\;A\in\{\alpha,i,{(\alpha)\atop (i)}\}$,
where ${\cal T}^B_A=G^{BD}{\cal T}_{DA}$, $"_{\vert A}"$, represents one from
the local covariant derivatives $"_{/\beta}"$, $"_{\vert j}"$ or
$"\vert^{(\beta)}_{(j)}"$, of the Cartan canonical connection $C\Gamma$.

In this context, let us denote
\begin{equation}
\begin{array}{lll}\medskip
\tilde{\cal T}_T=h^{\alpha\beta}\tilde{\cal T}_{\alpha\beta},&
\tilde{\cal T}_M=g^{ij}\tilde{\cal T}_{ij},&
\tilde{\cal T}_v=h_{\mu\nu}g^{mr}\tilde{\cal
T}^{(\mu)(\nu)}_{(m)(r)},\\
\tilde{\cal T}^\alpha_\beta=h^{\alpha\mu}\tilde{\cal T}_{\mu\beta},&
\tilde{\cal T}^i_j=g^{im}\tilde{\cal T}_{mj},&
\tilde{\cal T}^{(i)(\beta)}_{(\alpha)j}=h_{\alpha\mu}g^{mi}
\tilde{\cal T}^{(\mu)(\beta)}_{(m)(i)}.
\end{array}
\end{equation}

Following again the development of gravitational generalized metrical
multi-time theory from \cite{12}, we find
\begin{th}
If $p>2$, $n>2$, the new stress-energy d-tensors $\tilde{\cal T}_{\alpha\beta}$,
$\tilde{\cal T}_{ij}$ and $\tilde{\cal T}^{(\alpha)(\beta)}_{(i)(j)}$ of
$RGOGML^n_p$ must verify the following conservation laws:
\begin{equation}
\left\{\begin{array}{l}\medskip
\displaystyle{\tilde{\cal T}^\mu_{\beta/\mu}+{1\over 2-n}\tilde{\cal  T}_{M/\beta}
+{1\over 2-pn}\tilde{\cal T}_{v/\beta}=-R^m_{\beta\vert m}-P^{(m)}_{(\mu)\beta}
\vert^{(\mu)}_{(m)}}\\\medskip
\displaystyle{{1\over 2-p}\tilde{\cal T}_{T\vert j}+\tilde{\cal  T}^m_{j\vert m}
+{1\over 2-pn}\tilde{\cal T}_{v\vert j}=-P^{(m)}_{(\mu)j}\vert^{(\mu)}_{(m)}}\\
\displaystyle{{1\over 2-p}\tilde{\cal T}_T\vert^{(\alpha)}_{(i)}+{1\over 2-n}
\tilde{\cal  T}_M\vert^{(\alpha)}_{(i)}+\tilde{\cal T}^{(m)(\alpha)}_{(\mu)(i)}
\vert^{(\mu)}_{(m)}=-P^{m(\alpha)}_{\;\;\;\;i\vert m}},
\end{array}\right.
\end{equation}
\end{th}
where
\begin{equation}
\begin{array}{ll}\medskip
R^i_\beta=g^{im}R_{m\beta},&P^{(i)}_{(\alpha)\beta}=g^{im}h_{\alpha\mu}
P^{(\mu)}_{(m)\beta},\\
P^{i(\beta)}_{\;(j)}=g^{im}P_{m(j)}^{\;(\beta)},&
P^i_{(\alpha)j}=g^{im}h_{\alpha\mu}P^{(\mu)}_{(m)j}.
\end{array}
\end{equation}

\section{Maxwell equations of electromagnetic field}

\setcounter{equation}{0}
\hspace{5mm} In order to develope the electromagnetic theory on the generalized
metrical multi-time Lagrange space $RGOGML^n_p$, let us consider the {\it
canonical Liouville d-tensor} {\bf C}=$\displaystyle{x^i_\alpha{\partial\over
\partial  x^i_\alpha}}$
on $J^1(T,M)$. Using the Cartan canonical connection $C\Gamma$ of $RGOGML^n_p$,
we construct the {\it metrical deflection d-tensors} \cite{12}
\begin{equation}
\begin{array}{l}\medskip
\bar D^{(\alpha)}_{(i)\beta}=\left[G^{(\alpha)(\mu)}_{(i)(m)}x^m_\mu\right]_
{/\beta},\\\medskip
D^{(\alpha)}_{(i)j}=\left[G^{(\alpha)(\mu)}_{(i)(m)}x^m_\mu\right]_{\vert j},
\\\medskip
d^{(\alpha)(\beta)}_{(i)(j)}=\left[G^{(\alpha)(\mu)}_{(i)(m)}x^m_\mu\right]
\vert^{(\beta)}_{(j)},
\end{array}
\end{equation}
where $G^{(\alpha)(\beta)}_{(i)(j)}=h^{\alpha\beta}g_{ij}$
is the vertical fundamental  metrical d-tensor  of $RGOGML^n_p$ and
$"_{/\beta}"$, $"_{\vert j}"$ or $"\vert^{(\beta)}_{(j)}"$, are the local
covariant derivatives  of $C\Gamma$.
\medskip\\

Taking into account the expressions of the local covariant derivatives of
the Cartan canonical connection $C\Gamma$, we obtain
\begin{prop}
The metrical  deflection d-tensors of the space $RGOGML^n_p$ are given by
the following formulas:
\begin{equation}
\begin{array}{l}\medskip
\bar D^{(\alpha)}_{(i)\beta}=G^{(\alpha)(\mu)}_{(i)(m)}G^m_{r\beta}
x^r_\mu,\\\medskip
D^{(\alpha)}_{(i)j}=G^{(\alpha)(\mu)}_{(i)(m)}\Lambda^m_{rj}x^r_\mu,
\\\medskip
d^{(\alpha)(\beta)}_{(i)(j)}=G^{(\alpha)(\beta)}_{(i)(j)}+
C^\beta_{mji}h^{\alpha\mu}x^m_\mu.
\end{array}
\end{equation}
\end{prop}
\addtocounter{defin}{1}
{\bf Definition \thedefin} The distinguished 2-form on $J^1(T,M)$,
\begin{equation}
F=F^{(\alpha)}_{(i)j}\delta x^i_\alpha\wedge dx^j+f^{(\alpha)(\beta)}_{(i)(j)}
\delta x^i_\alpha\wedge\delta x^j_\beta,
\end{equation}
where
$F^{(\alpha)}_{(i)j}=\displaystyle{{1\over 2}\left[D^{(\alpha)}
_{(i)j}-D^{(\alpha)}_{(j)i}\right]}$ and
$f^{(\alpha)(\beta)}_{(i)(j)}=\displaystyle{{1\over 2}\left[d^{(\alpha)(\beta)}
_{(i)(j)}-d^{(\alpha)(\beta)}_{(j)(i)}\right]}$, is called the distinguished
{\it electromagnetic 2-form} of the generalized  metrical multi-time Lagrange
space $RGOGML^n_p$.
\begin{prop}
The local electromagnetic d-tensors of $RGOGML^n_p$ have the expressions,
\begin{equation}
\left\{\begin{array}{l}\medskip
\displaystyle{F^{(\alpha)}_{(i)j}=\left[\varphi_{ir}\Lambda^r_{mj}
-\varphi_{jr}\Lambda^r_{mi}+A_i\Lambda^0_{mj}-A_j\Lambda^0_{mi}\right]
h^{\alpha\mu}x^m_\mu},\\
\displaystyle{f^{(\alpha)(\beta)}_{(i)(j)}={1\over 2}\left[C^\beta_{mji}-
C^\beta_{mij}\right]h^{\alpha\mu}x^m_\mu,}
\end{array}\right.
\end{equation}
where $\Lambda^0_{mj}=\Lambda^r_{mj}A_r$.
\end{prop}

Particularizing the Maxwell equations of the electromagnetic field, described
in the general case of a generalized metrical multi-time Lagrange space
\cite{12}, we  deduce the main result of the electromagnetism on $RGOGML^n_p$.

\begin{th}
The electromagnetic components $F^{(\alpha)}_{(i)j}$ and $f^{(\alpha)(\beta)}_
{(i)(j)}$ of the generalized metrical multi-time Lagrange space $RGOGML^n_p$ are
governed by the Maxwell equations:
\medskip
\begin{equation}\hspace*{8mm}
\left\{\begin{array}{l}\medskip
\displaystyle{F^{(\alpha)}_{(i)k/\beta}={1\over 2}{\cal A}_
{\{i,k\}}\left\{\left[\bar D^{(\alpha)}_{(i)\beta}+x^{(\alpha)}_{(p)}
G^p_{i\beta}\right]_{\vert  k}-D^{(\alpha)}_{(i)m}G^m_{k\beta}
\right\}}\\\medskip
\displaystyle{f^{(\alpha)(\gamma)}_{(i)(k)/\beta}={1\over 2}{\cal A}_{\{i,k\}}
\left\{\bar D^{(\alpha)}_{(i)\beta}\vert^{(\gamma)}_{(k)}+
x^{(\alpha)}_{(p)}{\partial G^p_{i\beta}\over\partial x^k_\gamma}-\left[
d^{(\alpha)(\mu)}_{(i)(m)}+C^{p(\mu)}_{i(m)}x^{(\alpha)}_{(p)}\right]
G^m_{k\mu}\right\}}\\\medskip
\displaystyle{\sum_{\{i,j,k\}}F^{(\alpha)}_{(i)j\vert k}=-{1\over 2}\sum_
{\{i,j,k\}}}\left[d^{(\alpha)(\mu)}_{(i)(m)}+C^{p(\mu)}_{i(m)}x^{(\alpha)}_
{(p)}\right]r^m_{jks}x^s_\mu\\\medskip
\displaystyle{\sum_{\{i,j,k\}}\left\{F^{(\alpha)}_{(i)j}\vert^{(\gamma)}_{(k)}+
f^{(\alpha)(\gamma)}_{(i)(j)\vert k}\right\}=0}
\\
\displaystyle{\sum_{\{i,j,k\}}f^{(\alpha)(\beta)}_{(i)(j)}\vert^{(\gamma)}_
{(k)}=0,}
\end{array}\right.
\end{equation}
where $x^{(\alpha)}_{(p)}=G^{(\alpha)(\mu)}_{(p)(m)}x^m_\mu$.
\end{th}

\begin{center}
University POLITEHNICA of Bucharest\\
Department of Mathematics I\\
Splaiul Independentei 313\\
77206 Bucharest, Romania\\
e-mail: mircea@mathem.pub.ro\\
\end{center}

\end{document}